\documentclass[a4paper,10pt]{article}
\usepackage[cp1251]{inputenc}
\usepackage[english]{babel}

\usepackage{amssymb}\usepackage{amsmath}
\usepackage{epsfig}
\usepackage{bm}
\usepackage{graphics}
\def\bl{\rule[-1mm]{2.4mm}{2.4mm}}

\def\be{\begin{equation}}
\def\ee{\end{equation}}

\newtheorem{thrm}{\bf Theorem}

\newtheorem{lmm}{\bf Lemma}
\newtheorem{rmk}{\bf Remark}

\begin{document}

\title {Blaschke product for bordered surface}
\author{\copyright 2018 ~~~~A.B.Bogatyrev
\thanks{Supported by the Program of the Russ. Acad. Sci. 
'Fundamental Mathematics and its Applications' under grant PRAS-18-01}}
\date{}
\maketitle

\hfill {\it To the memory of Sasha Vassiliev}

It is well known that a ramified holomorphic covering of a closed unitary disc by another such a disc is given by a finite Blaschke product.
The inverse is also true. In this note we give two explicit constructions for a holomorphic ramified covering of a disc by other 
bordered Riemann surface. The problem of covering a disc by an annulus arises \cite{Kru} e.g. in multidimensional complex analysis;
we show that it may be effectively solved in terms of elliptic theta functions. The covering of a disc by a multiconnected flat domain is discussed in \cite{Gol}, Chap VI. The machinery used here strongly resembles the description of magnetic configurations in submicron planar magnets \cite{B17}.

\section{Schottky double, homologies and real differentials} 
Preliminary we fix some notations.

Let $X$ be a genus $g$ Riemann surface with $k>0$ boundary components called (real) ovals. It's Schottky double is 
a compact borderless surface $X_2$ of genus $g_2:=2g+k-1$ obtained from two copies of $X$ identified at the boundaries \cite{SS}. 
The anticonformal involution (reflection) $\tau$ naturally acts on the double interchanging the points on two copies of the bordered surface
with the ovals being the set of its fixed points.
  
\subsection{Cycles, including relative ones}
\begin{figure}
\centerline{\includegraphics[scale=.75, trim = 1cm 17cm 4.7cm 1cm, clip]{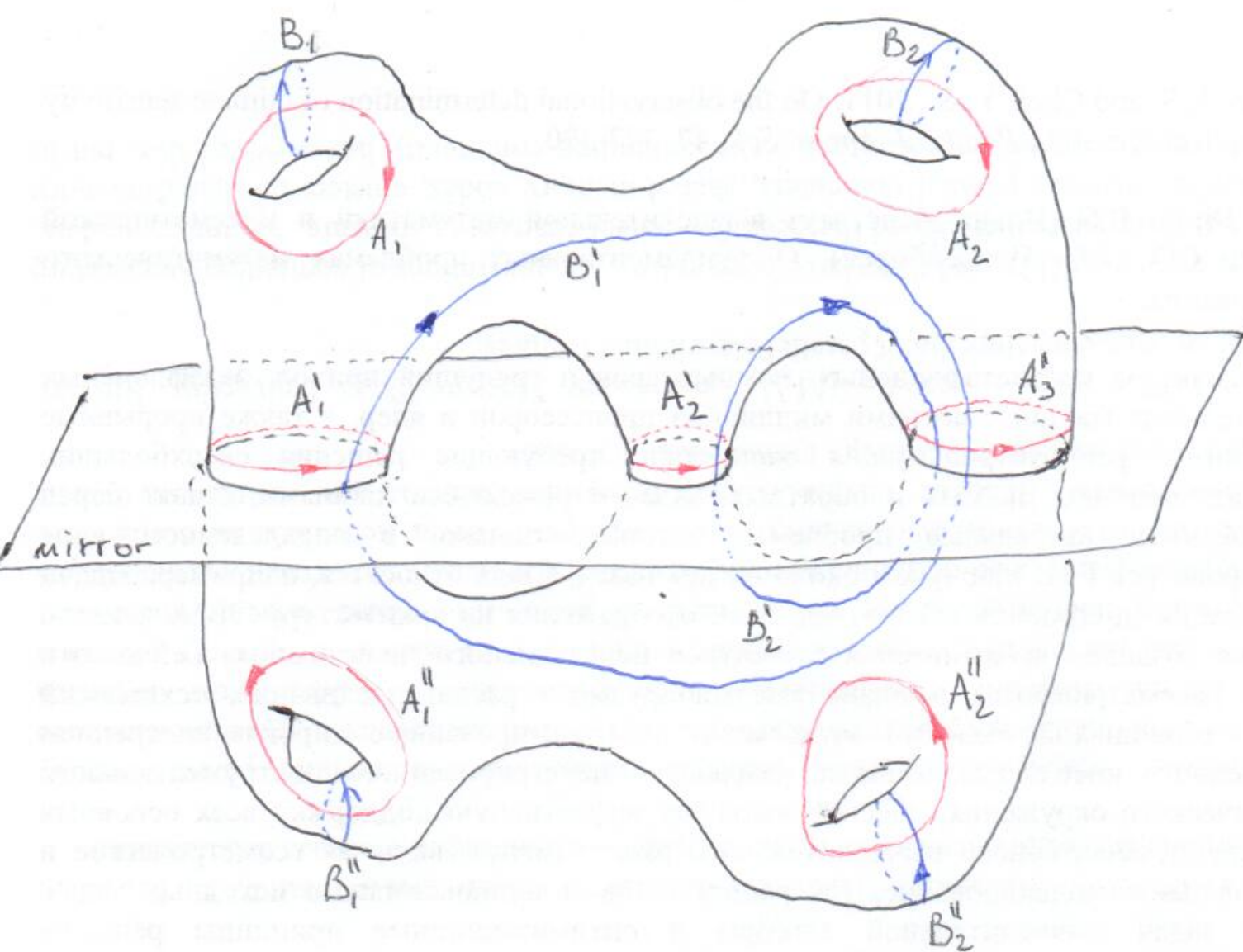}}
\caption{Fay homology basis on the double $X_2$ of the bordered surface $X$ with $g=2$ and $k=3$.} 
\label{CyclDbl}
\end{figure}

We fix a basis $A_1,\dots,A_g, B_1,\dots,B_g, A'_1,\dots,A'_{k-1} $ in the integer homology group $H_1(X,\mathbb{Z})$:
a standard pair of cycles $A_s,B_s$, $s=1,\dots,g$,  for each handle of the surface and all boundary components $A'_s$,
$s=1,\dots,k-1$, but one of them $A'_k$. The choice is by no means unique -- see Fig. \ref{CyclDbl} for the possible one.
Same list of cycles with removed boundary components $A'_j$ and amended by the arcs $B_j^+$, $j=1,\dots,k-1$, in $X$ which connect the oval $A'_j$ to  $A'_k$  and disjoint from other arcs of this kind as well as the basic cycles gives us a basis  
in the relative homology group $H_1(X,\partial X,\mathbb{Z})$. 

The lattice of boundary cycles $H_1(\partial X,\mathbb{Z})=\mathbb{Z}^k$ is naturally included to $H_1(X,\mathbb{Z})$, however with the loss of the rank: the sum of all boundary ovals with proper orientations is homological to zero. The integer cohomology lattice $H^1(X,\mathbb{Z})$ contains the sublattice $H_1^\perp(\partial X,\mathbb
{Z})$ of functionals vanishing on all boundary cycles of the surface.

The contours we introduced may be used to build  a symplectic basis in the homologies $H_1(X_2,\mathbb{Z})$ 
of the double (Cf: \cite{Fay}, Chap. 6):
\be
\label{FayHomo}
A''_s=\tau A_s, \quad B''_s=-\tau B_s,\quad s=1,\dots,g;
\qquad
B'_j=B_j^+-\tau B_j^+,\quad j=1,\dots,k-1.
\ee
Thus introduced basis possesses the standard intersection form and the easily checked behaviour under the reflection 
(Cf: Bobenko, Fay and Vinnikov homology basis):
\be
\label{FayTrans}
\begin{array}{l}
\tau A_s=A''_s;\qquad \tau B_s=-B''_s; \qquad s=1,\dots,g;\\
\tau A'_j=A'_j;\qquad \tau B'_j=-B'_j; \qquad j=1,\dots,k-1.
\end{array}
\ee

\subsection{Holomorphic differentials}
We consider the space $\mathbb{R}\Omega^1(X_2)\simeq\mathbb{R}^{g_2}$ of real holomorphic differentials on $X$, that is (holomorphic) 
differentials $d\rho$ on the double surface possessing the reflectional symmetry:
\be
\tau d\rho=\overline{d\rho}.
\ee
Riemann bilinear relations \cite{FK} guarantee that the following two pairings are non degenerate:
\be
\label{Recoupling}
\langle d\rho|C\rangle:= Re\int_C d\rho,   
\qquad d\rho\in\mathbb{R}\Omega^1(X_2),
\qquad C\in H_1(X,\mathbb{Z}),
\ee
\be
\label{Imcoupling}
\langle d\zeta|C\rangle:= Im\int_C d\zeta,   
\qquad d\zeta\in\mathbb{R}\Omega^1(X_2),
\qquad C\in H_1(X,\partial X,\mathbb{Z}).
\ee

Indeed, the symmetry of real differentials (possibly meromorphic) implies 
\be \label{ReflSym}
\int_{\tau C}d\rho=\overline{\int_C d\rho}, 
\ee
which in turn means that if all the values \eqref{Recoupling} (resp. values \eqref{Imcoupling}) 
for the given real differential $d\rho$ vanish then all periods of it on the double surface are purely imaginary (resp. real) and therefore $d\rho$ is zero. The last statement means that integer cohomologies $H^1(X,\mathbb{Z})$ and $H^1(X,\partial X, \mathbb{Z})$ may be naturally realized as rank $g_2$ lattice in the space of real differentials and in particular, there is a unique basis of differentials
dual to the above chosen basis in 1- (relative) homologies of $X$:
\be
\label{dualBasisRe}
Re~\int_C d\rho_D =\delta_{CD}, \qquad C,D\in \{A_1,\dots,A_g, B_1,\dots,B_g, A'_1,\dots,A'_{k-1}\},
\ee  
\be
\label{dualBasisIm}
Re~\int_C d\zeta_D =\delta_{CD}, \qquad C,D\in \{A_1,\dots,A_g, B_1,\dots,B_g, B^+_1,\dots,B^+_{k-1}\}.
\ee  
Differentials $d\rho_{A_s}, d\rho_{B_s}$, $s=1,\dots,g$ represent a basis for the mentioned above 
sublattice $H_1^\perp(\partial X,\mathbb{Z})$ of integer cohomologies.

\subsection{Meromorphic differentials}\label{MeroDiff}
On the double surface $X_2$ there exist a unique third kind abelian differential $d\eta_{zp}$
with simple poles at points $z,p$, residues $+1,1$ respectively and real normalization: all cyclic periods  of $d\eta_{zp}$ are purely imaginary. For the construction of meromorphic functions with mirror symmetry we use the differentials whose polar sets inherit the named  symmetry. We consider either (i) poles on real ovals of $X_2$: $z=\tau z$, $p=\tau p$ and it is easy to check that abelian differential $d\eta_{zp}$ will be real (meromorphic) in this case or (ii) the poles which are the reflections of one another $z=\tau p$ and the differential itself is now imaginary: $\tau\eta_{zp}=-\overline{\eta_{zp}}$. 
Those differentials may be expressed rather explicitly via the solution of the boundary value problem for a harmonic function in $X$, 
namely: (i) Neumann problem with two point-like sources on the boundary or (ii) the Dirichlet problem with one delta-source inside \cite{B17}.
Schottky-Klein prime form \cite{Fay, Mum} may be used for the efficient representation as well.  The following reciprocity laws are valid for the  differentials of this kind:  

\begin{lmm}\label{RBE}
The Riemann bilinear relations for the $Re-$normalized differential $d\eta_{zp}$ take the form\\
(i) Poles  $z=\tau z$, $p=\tau p$ are on real ovals: 

\be
\begin{array}{ll}
\int_{A_s}d\eta_{zp}=-\pi i Re~\int_p^z d\rho_{B_s},&\\ 
\int_{B_s}d\eta_{zp}=\pi i Re~\int_p^z d\rho_{A_s},& \qquad s=1,\dots,g;\\
\int_{A'_j}d\eta_{zp}=0,&\\ 
\int_{B'_j}d\eta_{zp}=2\pi i Re~\int_p^z d\rho_{A'_j},&\qquad j=1,\dots,k-1.
\end{array}
\ee
When poles $z,p$ lie on the same boundary oval the integrals in the r.h.s. of the equalities 
taken along the ovals are already real; in the third line the l.h.s integral 
is taken in the sense of principal value if necessary;
the integration path on the right does not cross  basic cycles $A_s,B_s,B'j$.

(ii) Poles $z=\tau p\in X$ are the reflections of each other:
\be
\begin{array}{ll}
\int_{A_s}d\eta_{zp}=-\pi i Im~\int_p^z d\zeta_{B_s},&\\ 
\int_{B_s}d\eta_{zp}=\pi i Im~\int_p^z d\zeta_{A_s},& \qquad s=1,\dots,g;\\
\int_{A'_j}d\eta_{zp}=-\pi i Im~\int_p^z d\zeta_{B^+_j},&\\ 
\int_{B'_j}d\eta_{zp}=0& \qquad j=1,\dots,k-1.
\end{array}
\ee
the integration path on the right is mirror-symmetric and does not cross chosen above basic cycles.
\end{lmm}

\begin{figure}
\centerline{\includegraphics[scale=.75,  trim = 1cm 20cm 7cm 1cm, clip]{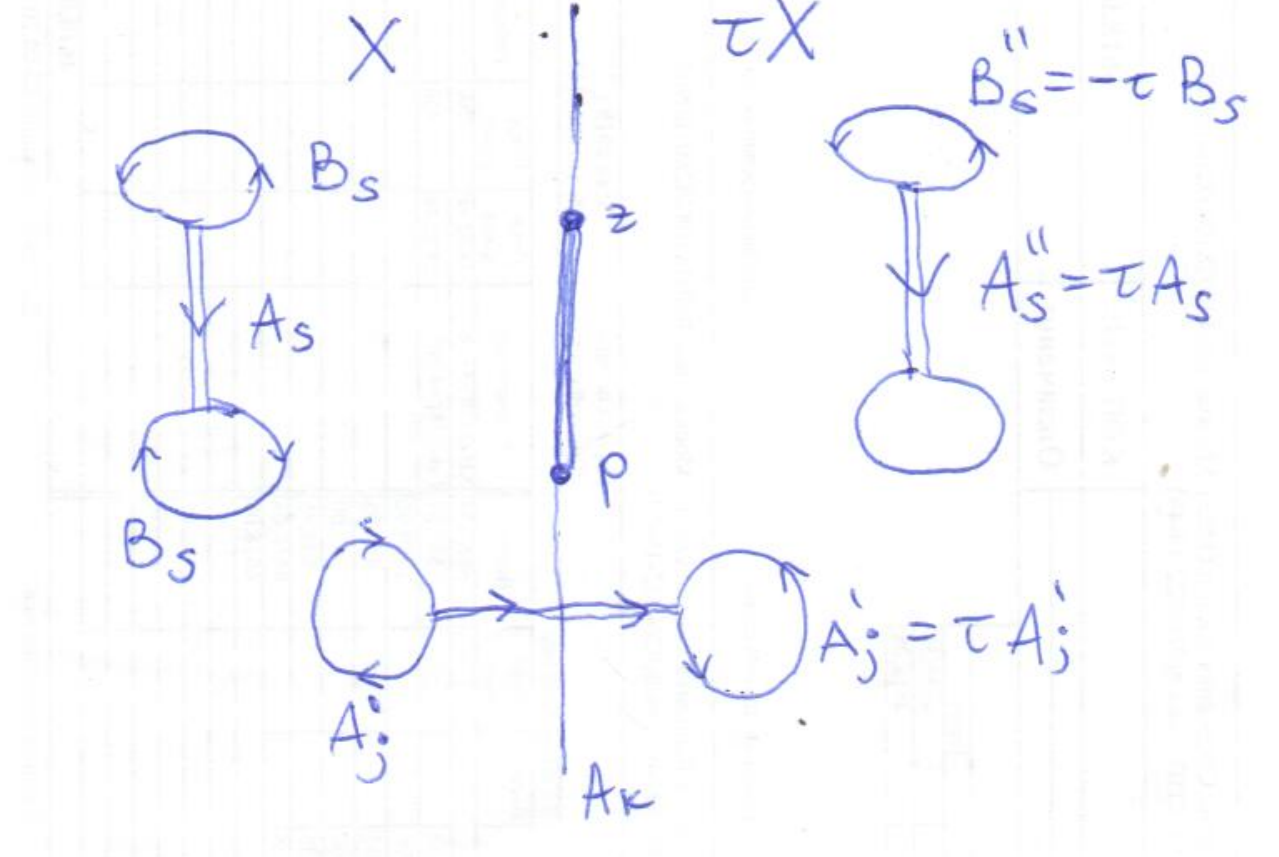}}
\caption{The surface $X_2$ cut along homology basis and an oval arc joining poles $z,p$.} 
\label{CutDbl}
\end{figure}

{\bf Proof} roughly follows the argumentation for the classical reciprocity laws \cite {FK, GH, Mum}
for holomorphic normalization. We give a sketch of it only.
We cut the surface $X_2$ along the chosen above basic cycles as well as the arc joining the poles $z,p$. In the arising flat surface one can introduce the single valued  abelian integral $\eta_{zp}$ -- see Fig. \ref{CutDbl}. Integrating  the product $\eta_{zp}d\rho$ with a real differential $d\rho$ from the basis \eqref{dualBasisRe} or respectively \eqref{dualBasisIm} along the boundary of the flat surface we get zero as the integrand has no singularities inside. On the other hand we can assemble the  boundary integrals into proper groups and get the above identities.

\section{Ramified coverings}
The main idea behind the description of ramified coverings $h$ of the closed upper half plane $\hat{\mathbb{H}}$ or the closed unit disc $\hat{\mathbb{D}}$ by a bordered surface $X$ is to extend the mapping to the double $X_2$ of the surface by reflection:
\be
\label{extH}
(i)\quad h(\tau x):=\overline{h(x)},   \quad ~~h(X)=\hat{\mathbb{H}},
\ee
\be
\label{extD}
(ii)\quad h(\tau x):=1/\overline{h(x)},   \quad h(X)=\hat{\mathbb{D}},
\ee
thus  getting a meromorphic function $h(x)\in\mathbb{C}(X_2)$. The zeros and poles of $h$ in the case (i) are simple, lie on the real ovals only and alternate on each oval since each boundary component of $X$ covers the extended real axis $\hat{\mathbb{R}}$. For the case (ii)
the divisor  of zeros and poles $(h):=(h)^+-(h)^-$, $(h)^\pm\ge0$, may have multiple points; all zeros $|(h)^+|$ lie in $Int~X$, while all poles $|(h)^-|$ lie the reflected surface: $(h)^-=\tau(h)^+$. The positions of points in the divisor $(h)$ are subjected to certain relations explained as follows: consider the differential 
\be
\label{deta}
d\eta_h:= d\log(h) = dh/h.
\ee
This differential is the third kind abelian differential with residue divisor $(h)$ and purely imaginary periods.
One easily checks that it is real in the case (i), imaginary in the case (ii) and may be decomposed to elementary 3rd kind abelian differentials introduced in Sect. \ref{MeroDiff}. All its periods lie in the lattice $2\pi i\mathbb{Z}$. The latter condition may be rewritten as restriction on the positions of points of $|(h)|$ using Riemann bilinear equations from Lemma \ref{RBE}. In each case (i), (ii)
the covering map admits the representation:
\be
\label{Coverformula}
h(u)=\exp(\int_v^u d\eta_h),
\ee
with $v\in h^{-1}(1)\subset\partial X$.

\begin{thrm} 
~({~\bf Covering of $\hat{\mathbb{H}}$})~
Let $h(u):~~X\to\hat{\mathbb{H}}$ be a ramified covering and $(h)=(h)^+-(h)^-$, $(h)^\pm\ge0$, be the decomposition of the divisor 
of its extension \eqref{extH} into zeros and poles, then the following restrictions hold:\\
(a) $|(h)^\pm|\subset \partial X$ and each oval contains at least one zero and one pole,\\
(b) Zeros and poles are simple and alternate on each oval,\\
(c) Their positions satisfy the lattice conditions:
\be
\label{latcondH}
\int_{(h)^-}^{(h)^+} d\rho\in \left\{
\begin{array}{ll}
\mathbb{Z},&\qquad d\rho\in H^1(X,\mathbb{Z})\subset\mathbb{R}\Omega^1(X_2),\\ 
2\mathbb{Z},&\qquad d\rho\in H_1^\perp(\partial X,\mathbb{Z})\subset H^1(X,\mathbb{Z}), 
\end{array}
\right.
\ee
here the cohomologies are realized as real differentials with the coupling \eqref{Recoupling}.

The converse is also true: once $(h)^+$ and $(h)^-$ are the divisors satisfying the above restrictions (a)-(c),
then formula \eqref{Coverformula} with $d\eta_h$ reconstructed from the divisor $(h):=(h)^+-(h)^-$ as above represents a ramified covering $X\to\hat{\mathbb{H}}$ provided $d\eta_h(v)>0$ (ovals are oriented as the boundary of $X$).
\end{thrm} 

\begin{thrm} 
~~~({\bf Covering of $\hat{\mathbb{D}}$ })~
Let $h(u):~~X\to\hat{\mathbb{D}}$ be a ramified covering and $(h)=(h)^+-(h)^-$, $(h)^\pm\ge0$, be the decomposition of the divisor 
of its extension \eqref{extD} into zeros and poles, then the following restrictions hold:\\
(a) $(h)^+\subset Int~X$, $(h)^+=\tau (h)^-$\\
(b) Their positions satisfy the lattice conditions:
\be
\label{latcondD}
Im~\int_{(h)^-}^{(h)^+} d\zeta\in 2\mathbb{Z},\qquad d\zeta\in H^1(X,\partial X,\mathbb{Z})\subset\mathbb{R}\Omega^1(X_2), 
\ee
here the cohomologies are realized as real  differentials with the pairing \eqref{Imcoupling}.

The converse is also true: once $(h)^+$ is the divisor satisfying the above restrictions (a)-(b),
then formula \eqref{Coverformula} with $d\eta_h$ reconstructed from $(h)$ as above represents a ramified covering $X\to\hat{\mathbb{D}}$ \end{thrm} 

\begin{rmk}
Upper half-plane may be conformally mapped to the unit dics and back by a suitable linear-fractional map. Hence we have two types of formulas for the covering mapping: the first of them is more convenient for the study of the mapping near the boundary, the other -- in the interior points of the bordered surface.   
\end{rmk}
\begin{rmk}
For the degree $n$ covering to exist we have to satisfy $g_2$ lattice conditions which may exceed the number of available degrees of freedom.
Sets of small winding numbers for the boundary components of $X$ may be not realizable.
A natural question arises \cite{Ahl, Or}:  what are possible sets of winding numbers for the ramified coverings of a given surface.
\end{rmk}

{\bf Proof of theorem 1.}  Conditions a) and b) are obvious since the restriction of $h:~\partial X\to  \partial\mathbb{H}$ is again a covering map: each oval winds around the real equator of the Riemann sphere. Lattice conditions (c) follow from the reciprocity laws of Lemma \ref{RBE}(i):
$$
2\pi i\int_{(h)^-}^{(h)^+}d\rho=2\pi i\sum_{z,p}\int_p^z d\rho= \int_C d\eta_h=Arg~ h(u)|_C\in 2\pi i \mathbb{Z},
\qquad d\rho\in H^1(X,\mathbb{Z}),
$$
where $(h)=:\sum_{z,p}(z-p)$ and each pair of points $z,p$ in this decomposition of the divisor lie on the same oval; 
the cycle $C$ equals to $2\sum_{s=1}^g(\langle d\rho|A_s\rangle B_s-\langle d\rho|B_s\rangle A_s)
+\sum_{j=1}^{k-1}\langle d\rho|A'_j\rangle B'_j$. If $d\rho$ annihilates boundary cycles, then $C\in 2H_1(X,\mathbb{Z})$ and the conditions \eqref{latcondH} follow.

We prove the converse under the assumptions (a)--(c) about  zeros and poles. In that case the following three statements are true.

1) Function $h(u)$ given by formula \eqref{Coverformula} with $Re-$normalized third kind differential $d\eta_h$
whose residue divisor is $(h)^+-(h)^-$ is single valued on the surface $X_2$. This is equivalent to 
\be
\label{inc}
\int_{H_1(X_2,\mathbb{Z})} d\eta_{h}\subset 2\pi i \mathbb{Z}.
\ee
The latter inclusion may be checked on the basic homology cycles. Taking into account the symmetry relation for real differentials 
\eqref{ReflSym} and the Fay basis transformation under the reflection \eqref{FayTrans}, it is enough to check the inclusion \eqref{inc} for the cycles $A_s$, $B_s$ and $B'_j$.
For $C=A_s, B_s$ from the reciprocity laws of Lemma 1 and the assumption (c) it follows:
$$
\int_C d\eta_h=\sum_{z,p} \int_Cd\eta_{zp}=
\pm \pi i \sum_{z,p} \int_z^p d\rho\in 2\pi i \mathbb{Z}, 
\qquad d\rho\in H_1^\perp(\partial X,\mathbb{Z}).
$$
For $C=B'_j$  we analogously have
$$
\int_C d\eta_h=\sum_{z,p} \int_Cd\eta_{zp}=
2\pi i \sum_{z,p} \int_z^p d\rho\in 2\pi i \mathbb{Z}, 
\qquad d\rho\in H^1(X,\mathbb{Z}).
$$

2) Function $h(u)$ is real since such is the differential $d\eta_h$ and $\tau v=v$: 
$$
h(\tau u):=\exp(\int_v^{\tau u} d\eta_h)=\exp(\int_v^u \tau d\eta_h)=\exp(\overline{\int_v^u d\eta_h})=\overline{h(u)}.
$$
In particular, $h(u)$ maps the boundary of $X$ to the boundary of $\mathbb{H}$.
 
3) $h^{-1}(\hat{\mathbb{R}})=\partial X$. Indeed, the points $h^{-1}\{0,\infty\}$ decompose the ovals into $2\deg h$
closed segments with disjoint interiors. The values of $h(u)$ in the interior of each segment have the same sign, for otherwise additional  zero or a pole of $h$ appears. Since $h$ has simple zeros and poles only, the neighbouring segments on ovals are mapped 1-1 to different segments $\pm[0,\infty]$ of the extended real line. Each real value has  at least $\deg h$ preimages in the ovals and hence no preimages in the interior of  $X$ or its reflection. 

Since $h$ presents a ramified covering of the Riemann sphere by the double surface $X_2$, its restriction to 
$X$ will be a cover of $\pm\hat{\mathbb{H}}$. To change the sign in the target space we just move the initial point $v$ to the neighbouring 
segment(s) tiling the ovals. ~~~\bl

{\bf Proof of Theorem 2}  follows the above scheme.  We only have to check that the image of $X$ under the mapping \eqref{Coverformula}
lies in the unit disc. To this end consider a function $W(u):=Re~\int_v^ud\eta_h, \qquad v\in \partial X,$ 
then 
$$W(\tau u)=Re~\int_v^{\tau u}d\eta_h=Re~\int_{\tau v}^{\tau u}d\eta_h=Re~\int_v^u\overline{d\eta_h}=-W(u).$$
Hence $W(\partial X)=0$. By maximum principle, $W(X)=[-\infty, 0]$ since $|(h)^+|\subset X$. So $h(X)\subset\mathbb{D}$.  ~~~\bl

\section{Examples}
\label {Examples}
\subsection{Disc} The disc has topological invariants $g=0,~k=1$ and the Theorem 1 above gives us a real rational function $R(u)$ with
alternating real zeros of numerator and denominator. This function may be transformed to a classical Blaschke product
by making a composition  $l\circ R\circ l^{-1}$ with linear-fractional function $l(u)$ mapping the upper half plane to the unit disc, 
say $l(u):=(u-i)/(u+i)$. Conjugation by linear fractional map sends $R$ to the same degree rational function with zeros 
at $l( R^{-1}(i))$ and poles at $l( R^{-1}(-i))$.  The latter two sets are mirror symmetric with respect to the unit circle  
(since $R(u)$ is real) and the first of them lies strictly inside the unit disc (since $Im~ R(u)$ is positive in the upper half plane only).
Hence we get a finite Blaschke product which Theorem 2 gives immediately. Note that the latter may have multiple zeros/poles unlike the representation given by Theorem 1.

\subsection{Annulus}
The annulus  $1\le|u|\le r$ has topological invariants $g=0$,  $k=2$. We can also represent it as a factor of the vertical strip $0\le Re(x)\le 1/2$ by a group of translations generated by $iT$, where $T=\pi/\log(r)>0$. The correspondence of the two models of the ring is given by the explicit formula $u(x)=\exp(2\pi x/T)$. The double of the annulus is the torus $X_2=\mathbb{C}/(\mathbb{Z}+iT\mathbb{Z})$
with the anticonformal involution $\tau x=-\bar{x}$.

\subsubsection{Cover of the half plane}
The generator of $H^1(X,\mathbb{Z})$ is $d\rho:=\frac{dx}{iT}$ and $H_1^\perp(\partial X,\mathbb{Z})$ is empty. 
Let the sets of zeros/poles be represented by the alternating points $x=z_j$ and $x=p_j$, $j=1,\dots,N$ on the sides on the strip. The Abel's lattice condition \eqref{latcondH} takes the form:  
\be 
\label{CondH}
\frac1{iT}\sum_{s=1}^N (z_s-p_s)\in\mathbb{Z}, 
\ee
which in terms of the concentric ring model exactly means that the product of all poles divided by the product of all zeros is positive. 
With this choice of representatives, the covering map $X\to\hat{\mathbb{H}}$ is proportional to the following
\be
h(x)=\exp(-2\pi imx)\prod_{j=1}^N\frac{\theta_1(x-z_j)}{\theta_1(x-p_j)},
\ee
here $\theta_1(x)=\exp(-\pi T/4)\sin(\pi x)-\dots$ is the only odd elliptic theta function of the modulus $iT$ \cite{Ach}
and $m$ is the integer number in the r.h.s. of the inclusion \eqref{CondH}. 

\subsubsection{Cover of the unit disc}
The generator of the relative cohomologies is $d\zeta=2idx$ and the Abel's lattice condition looks like 
\be
\label{CondD}
\sum_{s=1}^N (z_s-p_s)= 2Re~\sum_{s=1}^N z_s\in\mathbb{Z}, 
\ee
where $z_j$ represent zeros, lie in the mentioned above strip and $p_j=-\overline{z_j}$
are the poles. The $N-$ sheeted covering of a unit disc has the appearing 
\be
h(x)=\prod_{j=1}^N\frac{\theta_1(x-z_j)}{\theta_1(x+\overline{z_j})},
\ee
strongly resembling the classical Blaschke product. Now zeros and poles may be multiple.

In more sophisticated cases with $\chi(X)\le -1$, formula \eqref{Coverformula} for the ramified covering map can also be done computationally efficient with the use of higher genus Riemann theta functions \cite{Mum, Fay}. 


\vspace{5mm}
\parbox{8cm}
{\it
Marchuk Institute for Numerical Mathematics,\\
Russian Academy of Sciences;\\[1mm]

{\tt ab.bogatyrev@gmail.com}}

\end{document}